%% file: ConjClasses.tex
\documentclass{article}

\usepackage{latexsym,epsfig}
\usepackage{amssymb,amsmath,amsfonts}
\usepackage{exscale} 
\usepackage{ifthen}
\usepackage{longtable}
\usepackage{subfigure}
\usepackage{caption}
\usepackage{cite}
\usepackage{lscape}
\usepackage{alltt} 
\usepackage{multicol} 

\begin{document}
\protect\pagenumbering{arabic}
\setcounter{page}{1}
 
\input{commands.tex}

\title{\sc Asymptotic geometry and growth of conjugacy classes of nonpositively curved manifolds}
\author{\sc Gabriele Link}
\date{}
\maketitle
\begin{abstract} Let $\XX$ be a Hadamard manifold and
$\Gamma\subset\is(\XX)$ a discrete group of isometries which
contains an axial isometry without invariant flat half plane. We study the behavior of conformal densities on the
limit set of $\Gamma$ in order to derive a new asymptotic estimate for the
growth rate of closed geodesics in not necessarily compact or
finite volume manifolds.
\end{abstract}
\vspace{0.2cm}

\section{Introduction}

Let $M$ be a complete Riemannian manifold of nonpositive sectional
curvature, and denote by $P(t)$ the number of primitive closed geodesics in
$M$ of period $\le t$ modulo free homotopy. Our main interest in
this paper is the asymptotic behavior of this function.

In the case of a negative upper bound on the sectional curvature
of $M$, there is only one closed geodesic in each free homotopy
class. Using the ergodic theory of the geodesic flow,
G.~A.~Margulis (\cite{M},\cite{MS}) proved that for compact manifolds of
pinched negative curvature with volume entropy $h$
$$ P(t)h t e^{-ht} \ \to \ 1\qquad \mbox{as}\quad t\to\infty\,.$$
Recently, M.~Coornaert and G.~Knieper established an analogous
generalization of Margulis' result for compact quotients of
Gromov hyperbolic metric spaces (\cite[Theorem 1.1]{CK}).

G.~Knieper (\cite{K1},\cite{K2},\cite{K3}) extended the theory to compact geometric
rank one manifolds $M$.  If $h$ denotes the volume entropy of $M$,
and $P_{hyp}(t)$ the number of closed geodesics of period $\le t$
which do not admit a perpendicular parallel Jacobi field, he
proved the existence of constants $a>1$ and $t_0>0$ \st
$$ \frac1{a\;t}\;e^{ht}\le P_{hyp}(t)\le P(t)\le \frac{a}{t} e^{ht}
$$ for $t>t_0$ (\cite[Theorem~5.6.2]{K3}).

The purpose of this paper is a partial generalization of this result to a larger class of manifolds $M$. 
Let $\XX$ be the Riemannian universal covering
manifold of $M$, and $\Gamma$ the group of deck transformations.  Then $\XX$ is a Hadamard manifold, $M=\XX/\Gamma$, and $\Gamma$ is a discrete torsion free
subgroup of the isometry group $\is(\XX)$ of $\XX$. Let $\rand$ denote the geometric boundary of $\XX$ endowed with the
cone topology. We will only
require that $\Gamma$ contains an axial isometry which translates a geodesic without
flat half plane in $\XX$ (see \cite{B1} for precise definitions) and does not possess a global fixed  point in $\rand$. We emphasize that
we do not assume the manifold $M$ to be compact or of finite
volume. Instead of the volume entropy, we will therefore consider
the {\hl critical exponent} $\delta(\Gamma)$ of $\Gamma$ which is defined as the exponent of convergence of the Poincar{\'e} series $P^s(x,y)=\sum_{\ging} e^{-sd(x,\gamma y)}$. Furthermore,
if $\XX$ does not admit a quotient of finite volume, then the rank
rigidity result of Ballmann (\cite{B2}) and Burns-Spatzier
(\cite{BS}) does not require $M$ to be a geometric rank one
manifold, hence every geodesic in $\XX$ may bound a flat strip of finite
width. 

If $W$ is an open set in $\XX$, we denote   ${\cal P}(W)$ the set of free
  homotopy classes of closed geodesics in $M$ which possess a representative with a
  lift to $\XX$ intersecting $W$ nontrivially, and ${\cal P}_h(W)\subseteq {\cal
  P}(W)$ the subset of free homotopy classes containing a representative with
  a lift which intersects $W$ and does not bound a flat half plane. Notice that due to the
  possible occurrence of flat strips along each geodesic in $\XX$ there can be
  infinitely many closed geodesics in every free homotopy class of closed
  geodesics in $M$. However, all closed geodesics in the same free homotopy
  class have the same period. We will be interested in   the number of elements of period $\le t$ in ${\cal P}(W)$ and
  ${\cal P}_h(W)$ which we denote by $P(t;W)$
  and $P_h(t;W)$ respectively. Our main result is the following:
\begin{satz}
Let $W$ be a sufficiently large bounded open set in $\XX$, and
$\Gamma$ as above. Then
$$\delta(\Gamma)=\lim_{t\to\infty}\big(\frac1{t} \log P(t;W)\big)=
\lim_{t\to\infty}\big(\frac1{t} \log P_h(t;W)\big)\,.$$ 
If $\Gamma$ is ``weakly cocompact'' (see Definition~\ref{weakcoc}), then there
exist constants $b>1$, $R>0$ \st for $t>R$
$$ \frac1{b\;t}\;e^{\delta(\Gamma)t}\le P_h(t;W)\le P(t;W)\le b
e^{\delta(\Gamma)t}\,.$$ 
\end{satz}
Although we extend some of the methods from \cite{K1} and \cite{K2} 
to noncompact manifolds, our proof  of the lower bound in Theorem~1
avoids the use of Lemma~2.7 in \cite{K1}. Instead of that we make 
use of Theorem~2 and Corollary~\ref{twosets} below. 

Fix $\xo\in\XX$ and put
$N_\Gamma(R):=\#\{\gamma\in\Gamma\;|\,d(\xo,\gamo)<R\}$. 
A large part of the 
present paper is devoted to the study of the behavior of 
conformal densities on the geometric boundary $\rand$ and their
relation to $N_\Gamma(R)$.
Generalizing Lemma~4.4 in \cite{K2}, we
derive that every $\alpha$--dimensional conformal density
satisfies $\alpha\ge\delta(\Gamma)$. Using this fact and similar
arguments as 
T.~Roblin in \cite{R}, we are able to prove  
\begin{satz}
If $\Gamma$ is as above, 
 then $\,\lim_{R\to\infty}\big(\frac1{R}\log
N_\Gamma(R)\big)\,$ exists and equals $\,\delta(\Gamma)$.
\end{satz}
This theorem extends the main theorem in \cite{R} to 
manifolds of nonpositive curvature as above which are not necessarily
CAT$(-1)$. Furthermore, its  Corollary~\ref{twosets} will be
one of the key ingredients in our proof of the lower bound of Theorem~1.

The paper is organized as follows: In section~2 we recall some
basic facts about Hadamard manifolds and discrete groups of
isometries which contain an axial isometry without flat half plane.
In section~3 we introduce the concept of conformal densities and prove
a so-called shadow lemma, Theorem~\ref{shadowlemma}. This theorem gives an
idea of the local behavior of conformal densities and allows to
deduce asymptotic bounds on the exponential growth rate of the
number of orbit points of $\Gamma$. Section~4 is devoted to the proof
of Theorem~2 and its corollaries, which will be a key ingredient in
the proof of the lower bound of Theorem~1. 
In section~5, we investigate the asymptotic growth rate of
geometrically distinct closed geodesics modulo free homotopy
in a complete Riemannian manifold of nonpositive curvature
and prove Theorem~1.

\section{Axial isometries of Hadamard manifolds}

In this section we recall a few properties of Hadamard manifolds which possess
a geodesic without flat half plane. Most of the material can be found in (\cite{B1}).

Let $\XX$ be a complete simply connected Riemannian manifold of
nonpositive sectional curvature. The geometric boundary $\rand$ of
$\XX$ is the set of equivalence classes of asymptotic geodesic
rays endowed with the cone topology (see e.g. \cite[chapter~II]{B3}). This boundary is homeomorphic
to the unit tangent space of an arbitrary point in $\XX$, and
$\ganz:=\XX\cup \rand$ is homeomorphic to a closed $N$--ball in
$\RR^N$, where $N=\dim\XX$. Moreover, the isometry group of $\XX$
has a natural action by homeomorphisms on the geometric boundary.

All geodesics are assumed to have unit speed. For $x\in\XX$ and
$z\in\ganz$ we denote by $\sigma_{x,z}$ the unique unit speed
geodesic emanating from $x$ which contains $z$. Let
$D:=\{(x,x)\;|\, x\in\XX\}$ denote the diagonal in $\XX$. For
later use we introduce the continuous projection
\begin{equation}\label{pr}
\begin{array}{rcl} \pr:\quad \ganz\times\XX\setminus D &\to & \rand\\
(z,x)\quad &\mapsto &
\sigma_{x,z}(-\infty)\,.\end{array}
\end{equation} 
We say that two points $\xi$, $\eta\in\rand$ can be joined by a
geodesic, if there exists a geodesic $\sigma$ with extremities
$\sigma(-\infty)=\xi$ and $\sigma(\infty)=\eta$. A
geodesic $\sigma: \RR\to\XX$ is said to bound a {\hl flat strip} of width
$c\ge 0$ if there exists a totally geodesic isometric embedding $i:
[0,c]\times \RR \to\XX$ \st $i(0,t)=\sigma(t)$ for any $t\in\RR$.
\begin{df}\label{hypg}
A geodesic $\sigma$ in $\XX$ is called {\hd hyperbolic} if it does not
bound a flat strip of infinite width. In this case we call
$$c(\sigma):=\sup\{ c\ge 0\;|\, \sigma \mbox{ bounds a flat strip of width } c\}$$ the {\hd width} of the hyperbolic geodesic $\sigma$. 
\end{df}
Notice that the definition of hyperbolic geodesics given in \cite{K2} is more restrictive than Definition~\ref{hypg} above, since we allow that hyperbolic geodesics bound a flat strip of finite width. We refer to geodesics satisfying Knieper's hyperbolicity condition as rank one geodesics and remark that the width of a rank one geodesic is zero. The following lemma is a direct consequence of Lemma~2.1 in \cite{B1} and its proof.
\begin{lem}\label{hypgeod} 
Let  $\sigma$ be a hyperbolic geodesic of width $c(\sigma)\ge 0$ with extremities
$\sigma(-\infty)$ and $\sigma(\infty)$. Then for all $\eps>0$ there exist 
neighborhoods $U$, $V\subset\rand$ of $\sigma(-\infty)$, $\sigma(\infty)$ with closures $\overline U$, $\overline V$ homeomorphic to closed balls and $\overline U\cap\overline V=\emptyset$ \st any pair
of points $(\xi,\eta)\in U\times V$ can be joined by a geodesic. 
Furthermore, if $\tilde \sigma$ is a geodesic with extremities in $U$ and $V$,
then $\tilde\sigma$ is hyperbolic and
$d(\tilde\sigma,\sigma(0))< c(\sigma)+\eps$. 
\end{lem} 
\begin{df} 
An isometry $\gamma\neq\id$ of $\XX$ is called {\hd axial}, if there exists a constant
$l=l(\gamma)>0$ and a geodesic $\sigma$ \st
$\gamma(\sigma(t))=\sigma(t+l)$ for all $t\in\RR$. We call
$l(\gamma)$ the {\hd translation length} of $\gamma$, and $\sigma$
an {\hd axis} of $\gamma$. The boundary point
$\gamma^+:=\sigma(\infty)$ is called the {\hd attractive fixed
point}, and $\gamma^-:=\sigma(-\infty)$ the {\hd repulsive fixed
point} of $\gamma$. We further put
$\Ax(\gamma):=\{ x\in\XX\;|\, d(x,\gamma x)=l(\gamma)\}$.
\end{df}
We remark that $\Ax(\gamma)$ consists of the union of parallel geodesics
translated by $\gamma$, and 
$\overline{\Ax(\gamma)}\cap\rand$ is exactly the set of fixed points of
$\gamma$.

The following kind of axial isometries will play a crucial role in the sequel.
\begin{df}\label{hypaxiso}
An axial isometry is called {\hd hyperbolic axial} if it possesses
a hyperbolic axis. The {\hd width} $c(\gamma)$ of a  hyperbolic axial isometry $\gamma$ is defined by
 $$c(\gamma):=\sup\{ d(x, \sigma_{y,\gamma^+})\;|\, x, y\in\Ax(\gamma)\}\,.$$
\end{df}
Notice that if $\gamma$ is
hyperbolic axial, then $\gamma^+$ and $\gamma^-$ are the only
fixed points of $\gamma$, and every axial isometry commuting with $\gamma$
possesses the same set of invariant geodesics as $\gamma$. Furthermore, if $\sigma$ is an axis of $\gamma$, then $c(\gamma)\le 2 c(\sigma)$. 

Let us  recall some further properties
of hyperbolic axial isometries stated in Theorem~2.2 of \cite{B1}.
\begin{lem}\label{hypdyn} Suppose $\gamma$ is a
hyperbolic axial isometry. Then
\begin{enumerate}
\item[(1)] every point $\xi\in\rand\setminus\{\gamma^+\}$ can be joined
to $\gamma^+$ by a geodesic, and all these geodesics are hyperbolic,
\item[(2)] given neighborhoods $U$ of $\gamma^+$ and $V$ of $\gamma^-$,
there exists $N_0\in\NN$ \st $\gamma^n(\ganz\setminus V)\subset U$ and
$\gamma^{-n}(\ganz\setminus U)\subset V$ for all $n\ge N_0$.
\end{enumerate}
\end{lem}

For a discrete subgroup $\Gamma\subset\is(\XX)$ the
geometric {\hl limit set} of $\Gamma$ is defined by
$\Lim:=\overline{\Gamma\at x}\cap\rand$, where $x\in \XX$ is arbitrary. We say that two points $\xi$,
$\eta\in\rand$ are {\hl dual} \wrt $\Gamma$ if for all
neighborhoods $U$ of $\xi$ and $V$ of $\eta$ there exists $\ging$
\st $\gamma(\ganz\setminus U)\subset V$ and
$\gamma^{-1}(\ganz\setminus V)\subset U$. In this case, both $\xi$
and $\eta$ belong to $\Lim$.\\

From here on we will  assume that $\Gamma\subset\is(\XX)$ is a discrete
group which contains a hyperbolic axial isometry $h$  and does not possess a global fixed  point in $\rand$. The following proposition recalls the
part of Theorem~2.8 in \cite{B1} which applies to our groups
$\Gamma$.
\begin{prp}
For every neighborhood $U$ of a limit point $\xi\in\Lim$ there
exists $\gamma\in\Gamma$ \st $\gamma h^+\in U\setminus\{\xi\}$.
Moreover, the closure of $\Gamma\at\xi$ equals
$\Lim$, and any two points of $\Lim$ are dual \wrt $\Gamma$. 
\end{prp}

For $x\in\XX$ and $r>0$ we denote by $B_x(r)$ the open ball of
radius $r$ centered at $x$. Our first result states that for
$\xi\in\rand$ the projection $\pr_\xi :=\pr(\xi,\cdot)$ of a
sufficiently large ball in $\XX$ contains an open set in $\rand$.

The main difficulty in generalizing the analogous statement Lemma~3.5
in \cite{K2} to our situation consists in the fact that every hyperbolic geodesic may bound a flat strip. We therefore have to uniformly bound  the width of such flat strips along the geodesics in question.

\begin{lem}\label{grossproj}
For each $x\in\XX$ there exists a constant $\,r>0$ \st for all
$\xi\in\rand$ $\pr_\xi(B_x(r))$ contains an open set
$U\subset\rand$ with $U\cap\Lim\neq\emptyset$.
\end{lem}
\prf\  Let $x\in\XX$ arbitrary, fix a point $y\in\Ax(h)$ and let $c(h)\ge 0$
denote the width of $h$ (see Definition~\ref{hypaxiso}). The idea is to
construct a covering of $\rand$ by open sets.

Let $\eps>0$ and $U^\pm$ be disjoint neighborhoods of $h^\pm$ as in
Lemma~\ref{hypgeod}. That is, any two points in $U^+$, $U^-$ can
be joined by a geodesic, every such geodesic $\sigma$ is
hyperbolic and $d(y,\sigma)< c(h)+\eps$. For each
$\eta\in\rand\setminus (U^+\cup U^-)$ we denote by $\sigma_\eta$ a
hyperbolic geodesic connecting $\eta$ and $h^+$, and by $c(\sigma_\eta)$ the width of $\sigma_\eta$.  Let $V_\eta$ be
a neighborhood of $\eta$, $U_\eta\subseteq U^+$ a neighborhood of
$h^+$  \st
any two points in $U_\eta$, $V_\eta$ can be joined by a geodesic,
every such geodesic $\sigma$ is hyperbolic and
$d(\sigma_\eta(0),\sigma)< c(\sigma_\eta)+\eps$. Then $$\rand\subseteq
U^+\cup U^-\cup \bigcup_{\eta\in\rand\setminus(U^+\cup U^-)}
V_\eta \,,$$ 
and, since $\rand$ is compact, there exist finitely many points
$\eta_1,\eta_2,\ldots,\eta_n\in\rand$ \st 
$$ \rand\subseteq U^+\cup U^-\cup\bigcup_{i=1}^n V_{\eta_i} \,.$$
If $r:=c(h)+\eps+d(x,y)+\max_{1\le i\le n}\left( c(\sigma_{\eta_i}) +d(x,\sigma_{\eta_i}(0)\right)$, 
 then for any $\xi\in\rand$, 
the projection
$\pr_\xi(B_x(r))$ contains an open set in $\rand$: If $\xi\in
V_{\eta_i}$ for some $i\in\{1,2,\ldots n\}$, then $U_{\eta_i}$ is
the desired set, and $h^+\in U_{\eta_i}$ implies $U_{\eta_i}\cap
\Lim\neq\emptyset$. If $\xi\in U^-$, then $h^+\in U^+\subseteq
\pr_\xi(B_x(r))$, if $\xi\in U^+$, then $h^-\in U^-\subseteq
\pr_\xi(B_x(r))$.\qed\\

The following lemma states that the limit set can be covered by
finitely many $\Gamma$--trans\-lates of an appropriate open set in
$\rand$. 
\begin{lem}\label{gammatrans}
For any open subset $A\subset\rand$ with $A\cap\Lim\neq \emptyset$
there exists a finite set\\ 
\mbox{$\{\gamma_1,\gamma_2,\ldots,
\gamma_m\}\subset\Gamma$} \st
$$\Lim\subseteq \bigcup_{i=1}^m \gamma_i A\,.$$
\end{lem}
\prf\ Fix $\xi\in A\cap\Lim$. If $\xi=h^-$, we pick a neighborhood $U$ of $\xi$ \st $U\subseteq
A$ and $h^+\notin \overline{U}$. If $\Lim\subseteq \overline{U}\cup \{h^+\}$,
we choose $\gamma\in\Gamma$ \st $\gamma h^+\in U\setminus\{\xi\}$. Then
$h^+\in \gamma^{-1}U$, hence $\Lim\subseteq U\cup\gamma^{-1}U$. If
$\Lim\nsubseteq \overline{U}\cup \{h^+\}$, there exists $\eta\in
\Lim\setminus \overline{U}$, $\eta\neq h^+$. Since $\eta$ and $\xi$ are dual
\wrt $\Gamma$, for any neighborhood $V$ of $\eta$ there exists
$\gamma\in\Gamma$ \st $\gamma^{-1}(\ganz\setminus V)\subset U$ and
$\gamma(\ganz\setminus U)\subset V$, in particular $\ganz\setminus
V\subset\gamma U $. Choose neighborhoods $V$ of $\eta$, $W$ of $h^+$
sufficiently small so that the closures of the sets $U, V,W$ are
pairwise disjoint. By Lemma~\ref{hypdyn} (2) there exists
$n\in\NN$ \st $h^n(\ganz\setminus U)\subset W$ and
$h^{-n}(\ganz\setminus W)\subset U$, in particular $\ganz\setminus
W\subset h^n U$. We conclude
$$ \Lim\subseteq \ganz\setminus V \cup \ganz\setminus W \subset
\gamma U \cup h^n U\subseteq \gamma A\cup h^n A\,.$$ The case
$\xi=h^+$ is analogous.

If $\xi\notin\{h^+,h^-\}$, we choose a neighborhood $U\subseteq A$
of $\xi$ \st $h^+, h^-\notin \overline{U}$. Pick $\ging$ \st $\gamma h^+\in
U\setminus\{\xi\}$. Then $\gamma h\gamma^{-1}$ is hyperbolic
axial, and by the discreteness of $\Gamma$ (see \cite[Lemma~2.9]{B1}), 
$\gamma h^-\neq h^+$. Hence there exist neighborhoods
$V$ of $\gamma h^-$ and $W$ of $h^+$ \st the closures of $U,V,W$
are pairwise disjoint. Since $h^+$ and $\gamma h^+$ are dual,
there exists $g\in\Gamma$ \st $g(\ganz\setminus W)\subset U$.
Furthermore, there exists $n\in\NN$ \st $(\gamma
h\gamma^{-1})^n(\ganz\setminus V)\subset U$. We conclude

$\Lim\subseteq \ganz\setminus W\cup\ganz\setminus V\subset g^{-1}U\cup
(\gamma h\gamma^{-1})^{-n}U\subseteq g^{-1} A\cup \gamma
h^{-n}\gamma^{-1} A\,.$\qed

\section{Conformal densities}

Let $\XX$ be a Hadamard manifold with Riemannian distance $d$, $\xo\in\XX$ a
fixed base point, and $\Gamma\subset\is(\XX)$ a
discrete infinite subgroup. For $x,y\in\XX$, $s\in\RR$ we denote
by
$$  P^s(x,y)=\sum_{\gamma\in\Gamma} e^{-s d(x,\gamma y)} $$
the {\hl Poincar{\'e} series}. Its exponent of convergence
$\delta(\Gamma)$ is independent of $x,y$ by the triangle
inequality, and is called the {\hl critical exponent} of $\Gamma$. If
\begin{equation}\label{countdef}
N_\Gamma(R):=\#\{\gamma\in\Gamma\;|\, d(\xo,\gamo)<R\}\,,\qquad
\Delta N_\Gamma(R):=N_\Gamma(R)-N_\Gamma(R-1)\,,
\end{equation} 
then an easy calculation shows that
$$\delta(\Gamma)=\limsup_{R\to\infty}\big(\frac1{R}\log N_\Gamma(R)\big)=\limsup_{R\to\infty}\big(\frac1{R}\log \Delta N_\Gamma(R)\big)\,.$$
For $z\in \XX$ we consider the continuous map
\begin{equation}\label{Busdef}
\begin{array}{rcl} \bs_z\;:\quad \XX\times\XX&\to &\RR\\
(x,y)&\mapsto &d(x,z)-d(y,z)\,.\end{array}
\end{equation} 
This map extends continuously to the boundary via
$$ \bs_\eta(x, y)\,:= \lim_{s\to\infty}\big(d(x,\sigma(s))-d(y,\sigma(s))\big)\,,$$
where $\sigma$ is an arbitrary ray in the class of $\eta\in\rand$.
For $\xi\in\rand$, $y\in\XX$, the function \be \bs_\xi(\cdot ,
y):\quad \XX &\to & \RR\\
x &\mapsto & \bs_\xi(x, y)\ee is called the {\hl Busemann
function} centered at $\xi$ based at $y$.  It is independent of
the chosen ray $\sigma$ in the class of $\xi$ (see also \cite[chapter~II]{B3}).
\begin{df}\label{confdens}
Let $\MM^+(\rand)$ denote the cone of positive finite Borel
measures on $\rand$, and $\alpha>0$. An {\hd
$\alpha$--dimensional conformal density} is a continuous map
$$\ba{rcc}\mu\,:\quad\XX &\to &\MM^+(\rand)\\ x&\mapsto&\mu_x\ea$$
\vspace{0mm} with the properties
\be (i) && \supp(\mu_{o})\subseteq \Lim\,,\\
(ii) && \gamma*\mu_x=\mu_{\gamma^{-1} x}\ \ \mbox{for any}\
\gamma\in\Gamma\,,\ x\in\XX\,,\\
(iii)  &&\frac{d\mu_x}{d\mu_\xo}(\eta)=e^{\alpha
\bs_{\eta}(\xo,x)}\quad\mbox{for any}\ \,x\in\XX\,,\
\eta\in\supp(\mu_\xo)\,.\ee
\end{df}
The existence of a $\delta(\Gamma)$--dimensional conformal density
$\mu$ was proved by G.~Knieper (\cite[Lemma~2.2]{K2}) for
Hadamard manifolds and arbitrary discrete infinite isometry
groups. From his construction it follows that $\mu_x(\Lim)>0$ for
all $x\in\XX$. Our goal in this section is a generalization of
Lemma~4.4 in \cite{K2} which is valid for discrete groups with
compact quotient of geometric rank one. We will only require that
$\Gamma$ is a discrete isometry group which contains a hyperbolic
axial isometry and possesses infinitely many limit points.

Before we state our result, we present a few preliminary lemmata
needed in the proof. For the remainder of this section we will
assume that $\Gamma\subset\is(\XX)$ is a discrete group which
contains a hyperbolic axial isometry $h$ and does not globally fix a  point in $\rand$. Given Lemma~\ref{gammatrans}, the following lemma and its corollary are straightforward from Lemma~4.1 in \cite{K2}.
\begin{lem}\label{fullsupport}
Let $\mu$ be a conformal density and $x\in\XX$. Then
$\mu_x(\Lim)>0$ implies $\supp(\mu_x)=\Lim$.
\end{lem}
\prf\ Suppose $\xi\in\Lim$, $\xi\notin\supp(\mu_x)$. Let $U$ be a
neighborhood of $\xi$ \st $\mu_x(U)=0$. By Lemma~\ref{gammatrans}
there exists $\gamma_1,\gamma_2,\ldots,\gamma_m\in\Gamma$ \st
$\Lim\subseteq\bigcup_{i=1}^m \gamma_i U$. Hence
$$\mu_x(\Lim)\le \sum_{i=1}^m \mu_x(\gamma_i U)=\sum_{i=1}^m
\mu_{\gamma_i^{-1}x}(U)=0\,,$$ because $\mu_{\gamma_i^{-1}x}$,
$1\le i\le m$, is absolutely continuous \wrt $\mu_x$. \qed
\begin{cor}
If $\mu$ is a $\delta(\Gamma)$--dimensional conformal density, then
for any $x\in\XX$  $\supp(\mu_x)=\Lim$.
\end{cor}
For $x\in\XX$, $\xi\in\rand$ and $\eps>0$, we put $\
C_{x,\xi}^\eps:=\{z\in\ganz\;|\, \angle_x(z,\xi)<\eps\}$.

The following two lemmata are easy generalizations of Lemma~4.2 and
Proposition~3.6 of \cite{K2} to the noncompact case. 
\begin{lem}\label{grossmasse}
Fix $x\in\XX$ and let $\mu$ be a conformal density with
$\mu_x(\Lim)>0$. Then for any $\eps>0$ there exists a constant
$q=q(x,\eps)>0$ \st $\,\mu_x(C_{x,\xi}^\eps)>q\,$ for all
$\xi\in\Lim$.
\end{lem} \prf\ Suppose the contrary is true. Then there exists $\eps>0$ and a
sequence $(\xi_n)\subset\Lim$ \st
$\mu_x(C_{x,\xi_n}^{2\eps})\to 0$ as $n\to\infty$. Passing to a
subsequence if necessary, we may assume that $\xi_n$ converges to
a point $\xi\in\rand$.
Then there exists $N_0\in\NN$ \st $C_{x,\xi}^{\eps}\subseteq
C_{x,\xi_n}^{2\eps}$ for $n\ge N_0$, hence
$\mu_x(C_{x,\xi}^{\eps})=0$. Since the limit set is closed, we
further have $\xi\in\Lim$.
Arguing as in the proof of
Lemma~\ref{fullsupport} we obtain a contradiction to
$\mu_x(\Lim)>0$.\qed

\begin{lem} For $x\in\XX$ there exist constants $c_0>0$ and
$\eps>0$ \st for any $c>c_0$ and $y\in\XX\setminus B_x(c)$
$$ \pr_y(B_x(c))\supseteq C_{x,\xi}^\eps\cap\rand\quad \mbox{for
some}\ \xi\in\Lim\,.$$
\end{lem}
\prf\ Fix $x\in\XX$ and suppose the assertion is not true. Then
for all $n\in\NN$ there exists a point $y_n\in\XX\setminus B_x(n)$
\st $\pr_{y_n}(B_x(n))\nsupseteq C_{x,\xi}^{1/n}\cap\rand$ for all
$\xi\in\Lim$. 
Passing to a
subsequence if necessary, we assume that $(y_n)$ converges to a
point $\eta\in\rand$.

Then either for all $r>0$, $\eps>0$ and $\xi\in\Lim$
$C_{x,\xi}^\eps\cap\rand\nsubseteq \pr_\eta(B_x(r))$ in
contradiction to Lemma~\ref{grossproj}, or there exist constants
$r>0$, $\eps>0$ and a point $\xi\in\Lim$
\st $C_{x,\xi}^{2\eps}\cap\rand\subseteq \pr_\eta(B_x(r))$.
However, the continuity of the map $\pr: \ganz\times\XX\setminus
D\to\rand$ would then imply the existence of $N_0\in\NN$ \st
$C_{x,\xi}^\eps\subseteq \pr_{y_n}(B_x(r))$ for any $n>N_0$, in
contradiction to the choice of $y_n$ for $n>r$ and
$1/n<\eps$.\qed\\

We are finally able to prove the main theorem of this section.
\begin{thr}\label{shadowlemma}
Let $\alpha>0$ and $\mu$ an $\alpha$--dimensional
conformal density of positive and finite total mass. Then there exists a
constant $c_0>0$ with the following property: For $c>c_0$ there
exists a constant $D(c)>1$ \st for all $\gamma\in\Gamma$ with
$d(\xo,\gamo)>c$ we have
$$ \frac1{D(c)}\; e^{-\alpha d(\xo,\gamma \xo)}\le \mu_\xo(\pr_\xo(B_{\gamo}(c)))\le D(c) e^{-\alpha d(\xo,\gamma
\xo)}\,.$$ \end{thr} \prf\ Fix $c_0>0$ as in the previous lemma.
By Lemma~\ref{grossmasse} there exists a constant $q>0$ \st for
all $y\in\XX\setminus B_\xo(c)$ we have
$\mu_\xo(\pr_y(B_\xo(c)))>q$. Hence for any $\ging$ with
$d(\xo,\gamo)>c$
$$ q<\mu_\xo(\pr_{\gamma^{-1}\xo}(B_\xo(c)))<\mu_\xo(\rand)\,.$$
Furthermore, if $\eta\in\pr_{\xo}(B_{\gamo}(c)))$ then $0\le
d(\xo,\gamo)-\bs_{\eta}(\xo,\gamo)\le 2c$ by elementary geometric
estimates. For $\gamma\in\Gamma$ we abbreviate
$S_\gamma:=\pr_\xo(B_\gamo(c))$ and conclude $$
q<\mu_\xo(\gamma^{-1}
S_\gamma)=\mu_\gamo(S_\gamma)=\int_{S_\gamma}
d\mu_{\gamo}=\int_{S_\gamma} e^{\alpha
\bs_{\eta}(\xo,\gamo)}d\mu_\xo(\eta)\le  e^{\alpha
d(\xo,\gamo)}\mu_\xo(S_\gamma)\,.$$ Similarly we have
$$\mu_\xo(\rand)\ge \mu_\xo(\gamma^{-1}S_\gamma)=\int_{S_\gamma} e^{\alpha
\bs_{\eta}(\xo,\gamo)}d\mu_\xo(\eta)\ge e^{-2\alpha c} e^{\alpha
d(\xo,\gamo)}\mu_\xo(S_\gamma)\,$$ and summarize $\quad e^{-\alpha
d(\xo,\gamo)} q< \mu_\xo(S_\gamma)\le e^{-\alpha d(\xo,\gamo)}
e^{2\alpha c} \mu_\xo(\rand)\,.$\qed\\

The following proposition will be crucial 
in order to apply the methods developed by T.~Roblin in \cite{R}. 
\begin{prp}\label{dimconf}
If an $\alpha$--dimensional conformal density $\mu$ of positive and finite total
mass exists, then $\alpha\ge \delta(\Gamma)$.
\end{prp}
\prf\ Suppose $\mu$ is an $\alpha$--dimensional conformal density
of positive and finite total mass. Let $c>c_0$ with $c_0$ as in
Theorem~\ref{shadowlemma}, and $R>c$ arbitrary. Since $\Gamma$ is
discrete, every ball of radius $c+1$ in $\XX$ contains at most
$M=M(c)$ orbit points $\gamma\xo$. Hence every point in $\rand$ is
covered by at most $M$ sets $\pr_\xo(B_{\gamma\xo}(c))$, $R-1\le
d(\xo,\gamo)< R$, and therefore
$$\sum_{\begin{smallmatrix}\ging\\
R-1\le d(\xo,\gamo)<R\end{smallmatrix}}\mu_\xo(\pr_\xo(B_{\gamma\xo}(c)))\le M
\mu_\xo(\rand)\,.$$ Recall from~(\ref{countdef}) the definition of
$\Delta N_\Gamma(R)$. We conclude 
\be
\Delta N_\Gamma(R)\frac1{D(c)}\;e^{-\alpha R}&\le &\sum_{\begin{smallmatrix}\ging\\
R-1\le d(\xo,\gamo)<R\end{smallmatrix}}\frac1{D(c)}\;e^{-\alpha d(\xo,\gamo)}\\ 
&\le
&\sum_{\begin{smallmatrix}\ging\\
R-1\le d(\xo,\gamo)<R\end{smallmatrix}}
\mu_\xo(\pr_\xo(B_{\gamma\xo}(c)))\le M \mu_\xo(\rand)\,,\ee hence
$\quad \delta(\Gamma)=\limsup_{R\to\infty} \frac1{R}\log(\Delta
N_\Gamma(R))\le \alpha$.\qed
\begin{cor}\label{orbupcount}
There exists a constant $b>0$ \st $N_\Gamma(R)\le b
e^{\delta(\Gamma) R}$ for sufficiently large $R>0$. \end{cor}
\prf\ We compute as in the proof of the previous proposition
\be
\Delta N_\Gamma(R) \frac1{D(c)}\;e^{-\delta(\Gamma) R}&\le
&\sum_{\begin{smallmatrix}\ging\\
R-1\le d(\xo,\gamo)<R\end{smallmatrix}}\frac1{D(c)}\;e^{-\delta(\Gamma) d(\xo,\gamo)}\\
&\le &\sum_{\begin{smallmatrix}\ging\\
R-1\le d(\xo,\gamo)<R\end{smallmatrix}}
\mu_\xo(\pr_\xo(B_{\gamma\xo}(c)))\le M \mu_\xo(\rand)\,,\ee hence
$\Delta N_\Gamma(R)\le M D(c)\mu_\xo(\rand) e^{\delta(\Gamma)R}$.
Furthermore, if $n$ denotes the smallest integer greater than $R$,
we have 
\be N_\Gamma(R)&\le & \sum_{j=1}^{n}\Delta N_\Gamma(j)\le
M D(c)\mu_\xo(\rand)\sum_{j=1}^n (e^{\delta(\Gamma)})^j\\
&= & M D(c)\mu_\xo(\rand)\frac{e^{\delta(\Gamma)}(e^{\delta(\Gamma)n}-1)}{e^{\delta(\Gamma)}-1}
\ee 
and the assertion follows with $b=M D(c)\mu_\xo(\rand)
e^{2\delta(\Gamma)}/(e^{\delta(\Gamma)}-1)$.\qed\\

We next introduce a class of groups for which we will be able to
derive stronger results.
\begin{df}\label{weakcoc}
A discrete subgroup $\Gamma\subset\is(\XX)$ is called {\hd weakly cocompact} if there
exists a $\,\delta(\Gamma)$--dimensional conformal density $\mu$ and
constants $b>0$, $c_\Gamma>0$ and $r>0$ \st for all $c\ge c_\Gamma$
$$ \hspace{2cm} \liminf_{R\to\infty} \mu_o\big(\hspace{-.6cm}\bigcup_{\stackrel{\gamma\in\Gamma}{R-r\le
d(\xo,\gamo)<R}}\hspace{-.6cm} \pr_\xo(B_\gamo(c))\big)\ge
b\,.$$\end{df} Notice that for discrete isometry groups of
Hadamard manifolds, cocompact implies weakly cocompact,
because in this case
$$0<\mu_\xo(\rand)=\mu_\xo(\hspace{-.5cm}\bigcup_{\begin{smallmatrix}\ging\\
R-1\le d(\xo,\gamo)<R\end{smallmatrix}}\hspace{-.5cm} \pr_\xo(B_\gamo(c))\big)$$ if $c\ge
\diam(\XX/\Gamma)$ and $R>c$. Further examples of weakly cocompact
groups are convex cocompact isometry groups
of real hyperbolic spaces, and radially cocompact isometry groups
(see \cite{L} for a definition) of symmetric spaces.

For weakly cocompact groups, we have the following lower bound for $N_\Gamma(R)$.
\begin{lem}\label{orblowcount}
Let $\Gamma\subset\is(\XX)$ be a weakly cocompact discrete group which
contains a hyperbolic axial isometry. Then there exists $a>0$ and
$R_0>0$ \st for all $R\ge R_0$
$$ N_\Gamma(R)\ge a e^{\delta(\Gamma)R}\,.$$
\end{lem}
\prf\  Let $c>\max\{c_0,c_\Gamma\}$ with $c_0>0$ as in
Theorem~\ref{shadowlemma}, $b>0$, $c_\Gamma>0$ and $r>0$ the constants
from the definition. Then there exists $R_0>r$ \st for all $R\ge R_0$
\be b & \le &
\mu_\xo\big(\hspace{-.5cm}\bigcup_{\begin{smallmatrix}\ging\\
R-r\le d(\xo,\gamo)<R\end{smallmatrix}} \hspace{-.5cm}\pr_\xo(B_\gamo(c))\big)\le
\sum_{\begin{smallmatrix}\ging\\
R-r\le d(\xo,\gamo)<R\end{smallmatrix}}\mu_\xo(\pr_\xo(B_{\gamo}(c)))\\
&\le & D(c) \sum_{\begin{smallmatrix}\ging\\
R-r\le d(\xo,\gamo)<R\end{smallmatrix}} e^{-\delta(\Gamma)d(\xo,\gamo)}\le D(c) N_\Gamma(R)
e^{-\delta(\Gamma)(R-r)}\,,\ee hence $N_\Gamma(R)\ge
 a e^{\delta(\Gamma)R}\ $ for $\,a=b e^{-r\delta(\Gamma)}/D(c)$.\qed

\section{The critical exponent}

In this section we are going to prove Theorem~\ref{deltaequal} and some
corollaries, which will be necessary to derive the lower bound in  
Theorem~\ref{main}. As before, $\XX$ will denote a Hadamard manifold,
$\xo\in\XX$ a fixed base point, and $\Gamma\subset\is(\XX)$ a discrete group
which contains a hyperbolic axial isometry and possesses infinitely many limit
points. Given Proposition~\ref{dimconf}, the arguments of T.~Roblin
in the context of CAT$(-1)$--spaces (\cite{R})  remain valid in our
setting. We include the proofs  for the convenience of the reader.
\begin{thr}\label{deltaequal}
If $\Gamma$ is a discrete isometry group of a Hadamard
manifold $\XX$ which contains a hyperbolic axial isometry and
possesses infinitely many limit points, and $\xo\in\XX$ a fixed base point, then
$\, \lim_{R\to\infty}\big(\frac1{R}\log N_\Gamma(R)\big)$ exists and equals
$\delta(\Gamma)$.
\end{thr}
\prf\ Assume that $\liminf_{R\to\infty}\big(\frac1{R}\log
N_\Gamma(R)\big)<\delta(\Gamma)$. Then there exists
 a sequence $(R_k)\subset\RR$,
$R_k\to\infty$, and $0<\alpha<\delta(\Gamma)$ \st
$N_\Gamma(R_k)\le e^{\alpha R_k}$ for all $k\in\NN$. We are going
to construct an $\alpha$--dimensional conformal density in order to
obtain a contradiction to Proposition~\ref{dimconf}.

Let $\delta$ denote the unit Dirac point measure. For $R>0$ and
$x\in\XX$ we put
$$\nu_x^R:=\sum_{\begin{smallmatrix}\ging\\
d(x,\gamo)\le R\end{smallmatrix}}e^{-\alpha
d(x,\gamo)}\delta(\gamo)/\big(\sum_{\begin{smallmatrix}\ging\\
d(\xo,\gamo)\le R\end{smallmatrix}}e^{-\alpha d(o,\gamo)}\big)\,.$$ From $\Vert \nu_o^R\Vert=1$ and
the Theorem of Banach-Alaoglu it follows that for any $r>0$ there
exists a sequence $k_n(r)\subset\NN$, $k_n(r)\to\infty$, \st
$\nu_o^{R_{k_n(r)}-r}$ converges weakly to a probability measure
$\mu_o^r$. Furthermore, the support of $\mu_o^r$ equals $\Lim$,
because, since $\alpha<\delta(\Gamma)$, the series in the
denominator diverges as $R\to\infty$. We denote by $\mu^r$ the
$\alpha$--dimensional conformal density induced by $\mu_o^r$. Our
aim is to prove that for any $x\in B_o(r)$ the measures
$\nu_x^{R_{k_n(r)}-r}$ converge weakly to $\mu_x^r$. That is, we
have to show that for every bounded and continuous function $f$ on
$\ganz$
$$\lim_{n\to\infty} \int_{\ganz} f(z)d\nu_x^{R_{k_n(r)}-r}(z)
=\int_{\ganz} f(z)d\mu_x^r(z)\,.$$ 
From the definition of
$\mu_x^r$ and $\mu_\xo^r$ it follows that
$$\int_{\ganz} f(z)d\mu_x^r(z)=\int_{\ganz} f(z)e^{\alpha \bs_z(o,x)} d\mu_o^r(z)=\lim_{n\to\infty} \int_{\ganz}
f(z)e^{\alpha \bs_z(o,x)}d\nu_o^{R_{k_n(r)}-r}(z)$$ with
$\bs_z(o,x)=d(o,z)-d(x,z)$, $z\in\ganz$, as in~(\ref{Busdef}). 
We put $\Vert f\Vert:=\sup_{x\in\ganz}|f(x)|$,
$$\Psi_n:=\sum_{\begin{smallmatrix}\ging\\
d(\xo,\gamo)\le R_{k_n(r)}-r\end{smallmatrix}}
e^{-\alpha d(o,\gamo)}$$ 
and compute for all $n\in\NN$ 
\be && \Big\vert \int_{\ganz}f(z)d\nu_x^{R_{k_n(r)}-r}(z)-
\int_{\ganz} f(z)e^{\alpha
\bs_z(o,x)}d\nu_o^{R_{k_n(r)}-r}(z)\Big\vert \\
&\le &\frac{\Vert f\Vert}{\Psi_n}\ \Big\vert\hspace{-.5cm}
\sum_{\begin{smallmatrix}\ging\\
d(x,\gamo)\le R_{k_n(r)}-r\end{smallmatrix}} e^{-\alpha
d(x,\gamo)}-\sum_{\begin{smallmatrix}\ging\\
d(\xo,\gamo)\le R_{k_n(r)}-r\end{smallmatrix}}
e^{\alpha(d(o,\gamo)-d(x,\gamo))}e^{-\alpha
d(o,\gamo)}\Big\vert\,.\ee Consider the number
$$ g_x(R):=\Big\vert
\sum_{\begin{smallmatrix}\ging\\
d(x,\gamo)\le R\end{smallmatrix}} e^{-\alpha
d(x,\gamo)}-\sum_{\begin{smallmatrix}\ging\\
d(\xo,\gamo)\le R\end{smallmatrix}} e^{-\alpha
d(x,\gamo)}\Big\vert\,.$$
Then for $x\in B_o(r)$ we have \be g_x(R)& \le &
\sum_{\begin{smallmatrix}\ging\\
|d(o,\gamo)-R|\le d(o,x)\end{smallmatrix}} e^{-\alpha
d(x,\gamo)}\le \sum_{\begin{smallmatrix}\ging\\
|d(o,\gamo)-R|\le d(o,x)\end{smallmatrix}}
e^{-\alpha
d(o,\gamo)+\alpha d(o,x)}\\
&\le & \sum_{\begin{smallmatrix}\ging\\
|d(o,\gamo)-R|\le d(o,x)\end{smallmatrix}}
e^{-\alpha R+2\alpha d(o,x)} \le N_\Gamma(R+r) e^{-\alpha R+
2\alpha r}\,.\ee We conclude 
\be &&\Big\vert
\int_{\ganz}f(z)d\nu_x^{R_{k_n(r)}-r}(z)- \int_{\ganz}
f(z)e^{\alpha \bs_z(o,x)}d\nu_o^{R_{k_n(r)}-r}(z)\Big\vert \\&\le
&\frac{\Vert f\Vert}{\Psi_n}\ g_x(R_{k_n(r)}-r)\le \frac{\Vert
f\Vert}{\Psi_n}\, N_\Gamma(R_{k_n(r)}) e^{-\alpha
R_{k_n(r)}+3\alpha r}\le \frac{\Vert f\Vert}{\Psi_n}\, e^{3\alpha
r}\to 0 \,,\ee because $\Psi_n$ is unbounded as $n\to\infty$. This
implies that $\nu_x^{R_{k_n(r)}-r}$ converges weakly to $\mu_x^r$
for all $x\in B_o(r)$.

Obviously, we have $\gamma*\nu_x^R=\nu_{\gamma^{-1}x}^R$ for all
$x\in\XX$, $\gamma\in\Gamma$. If $d(o,x)<r$ and $d(o,\gamma x)<r$,
this implies $\gamma*\mu_x^r=\mu_{\gamma^{-1}x}^r$.

We finally consider a sequence $(r_j)\subset\RR$, $r_j\to\infty$,
\st $\mu_o^{r_j}$ converges weakly to a probability measure
$\mu_o$. Let $\mu$ be the conformal density induced by $\mu_o$.
Then $\mu_x^{r_j}$ converges weakly to $\mu_x$ for all $x\in\XX$.
Furthermore $\gamma*\mu_x=\mu_{\gamma^{-1}x}$ for all $x\in\XX$,
$\gamma\in\Gamma$. This yields the desired 
contradiction.\qed\\

If $A\subset\rand$, $z\in\ganz$, we let
$\angle_o(z,A):=\inf_{\eta\in A}\angle_o(z, \eta)$ and
$$N_\Gamma(R;A):=\#\{\gamma\in\Gamma\;|\, d(o,\gamo)<R,\;
\angle_o(\gamo,A)=0\}\,.$$
\begin{cor}
If $A\subset\rand$ is an open set with $A\cap\Lim\neq\emptyset$,
then 
$$\lim_{R\to\infty}\big(\frac1{R}\log
N_\Gamma(R;A)\big)=\delta(\Gamma)\,.$$ 
For weakly cocompact $\Gamma$
there exist $b>1$ and $R_0>0$ \st for all $R>R_0$
$$ \frac1{b}e^{\delta(\Gamma)R}\le N_\Gamma(R;A)\le b
e^{\delta(\Gamma)R}\,.$$
\end{cor}
\prf\ Let $U\subset\rand$ be an open set with
$U\cap\Lim\neq\emptyset$ \st $\overline{U}\subset A$. By
Lemma~\ref{gammatrans} there exist
$\gamma_1,\gamma_2,\ldots,\gamma_m\in\Gamma$ \st
$\Lim\subseteq\bigcup_{i=1}^m \gamma_i U$. Let
$M\subseteq\bigcup_{i=1}^m \gamma_i U$ be an open set which
contains $\Lim$. Then $\angle_o(\gamo,M)=0\,$ for all but finitely
many $\gamma\in\Gamma$ by the definition of the limit set, hence
$N_\Gamma(R; M)\ge N_\Gamma(R)-j$ for some constant $j\in\NN$.

Fix $g\in\{\gamma_1,\gamma_2,\ldots,\gamma_m\}$ and suppose $\angle_\xo(\gamo,gU)=0$ and $\angle_\xo(g^{-1}\gamo,A)>0$ for infinitely many
$\gamma\in\Gamma$. Let $(\gamma_k)\subset\Gamma$ be a sequence
with this property. Passing to a subsequence if necessary, we may
assume that $\gamma_k\xo$ converges to a point $\eta\in g
\overline U$. Since $g\overline U\subset g A$, this
implies $\angle_{g \xo}(\gamma_k\xo,g A)=0$ for
almost all $k\in\NN$ in contradiction to
$\angle_{g\xo}(\gamma_k\xo,g A)=\angle_\xo(
g^{-1}\gamma_k\xo,A)>0\,.$ Hence
$$c(g):=\#
\{\gamma\in\Gamma\;|\, \angle_{\xo}(\gamo,g U)=0\ \an\
\angle_o(g^{-1}\gamo,A)>0\}<\infty\,,$$ and
$N_\Gamma(R;g U)\le
N_\Gamma(R+d(\xo,g\xo);A)+c(g)\,.$ Put $r:=\max_{1\le i\le
m} d(\xo,\gamma_i\xo)$ and $c:=\sum_{i=1}^m c(\gamma_i)$. Then \be
N_\Gamma(R-r;A)-j & \le & N_\Gamma(R-r)-j  \le
N_\Gamma(R-r;M)\le\sum_{i=1}^m N_\Gamma(R-r;\gamma_i U)\\
&\le & \sum_{i=1}^m \big(
N_\Gamma(R-r+d(\xo,\gamma_i\xo);A)+c(\gamma_i)\big)\le
mN_\Gamma(R;A)+c\,,\ee which proves the assertion. The claim for
weakly cocompact $\Gamma$ follows with Corollary~\ref{orbupcount}
and Lemma~\ref{orblowcount}.\qed\\

The second corollary of Theorem~\ref{deltaequal} estimates the
numbers
$$N_\Gamma(R;A,B):=\#\{\gamma\in\Gamma\;|\, d(o,\gamo)<R,\;
\angle_o(\gamo,A)=0,\;\angle_o(\gamma^{-1}\xo,B)=0\}\,,$$ where
$A,B\subseteq\rand$ are open sets.
\begin{cor}\label{twosets}
If $ A, B\subset\rand$ are  open sets with
$A\cap\Lim\neq\emptyset$, $B\cap\Lim\neq\emptyset$, then
$$\lim_{R\to\infty}\big(\frac1{R}\log
N_\Gamma(R;A,B)\big)=\delta(\Gamma)\,.$$
For weakly cocompact $\Gamma$ there exist $b>1$ and $R_0>0$ \st
for all $R>R_0$
$$ \frac1{b}e^{\delta(\Gamma)R}\le N_\Gamma(R;A,B)\le b
e^{\delta(\Gamma)R}\,.$$
\end{cor}
\prf\ Let $V\subset\rand$ be an open set with
$V\cap\Lim\neq\emptyset$ \st $\overline{V}\subset B$. By
Lemma~\ref{gammatrans} there exist
$\gamma_1,\gamma_2,\ldots,\gamma_m\in\Gamma$ \st
$\Lim\subseteq\bigcup_{i=1}^m \gamma_i V$. As in the proof of the previous corollary we let $M\subseteq\bigcup_{i=1}^m\gamma_i V$ be an open set which contains $\Lim$. Then there exists an integer $j\in\NN$ \st $N_\Gamma(R;A,M)\ge N_\Gamma(R;A)-j$. 

Fix $g\in\{\gamma_1,\gamma_2,\ldots,\gamma_m\}$ and suppose there exist
infinitely many $\gamma\in \Gamma$ \st $\angle_\xo(\gamma^{-1}\xo,g V)=0$
and $\angle_\xo(g^{-1}\gamma^{-1}\xo,B)>0$. Denote by
$(\gamma_k)\subset\Gamma$ be a sequence with this property.
Passing to a subsequence if necessary, we may assume that
$\gamma_k^{-1}\xo$ converges to a point $\eta\in g \overline
V$. Since $g\overline V\subset g B$, this implies
$\angle_{g \xo}(\gamma_k^{-1}\xo,g B)=0$ for almost
all $k\in\NN$, in contradiction to
$\angle_{g\xo}(\gamma_k^{-1}\xo,g B)=\angle_\xo(
g^{-1}\gamma_k^{-1}\xo,
B)>0\,.$ Hence
$$c(g):=\#
\{\gamma\in\Gamma\;|\, \angle_{\xo}(\gamma^{-1}\xo,g V)=0\
\an\ \angle_\xo(g^{-1}\gamma^{-1}\xo,B)>0\}<\infty\,,$$ and
$N_\Gamma(R;A,g V)\le
N_\Gamma(R+d(\xo,g\xo);A,B)+c(g)\,.$ We put\\
 $r:=\max_{1\le i\le
m} d(\xo,\gamma_i\xo)$ and $c:=\sum_{i=1}^m c(\gamma_i)$ and conclude
\be N_\Gamma(R-r;A,B)-j & \le & N_\Gamma(R-r;A)-j\le N_\Gamma(R-r;A,M)\\
&\le &
\sum_{i=1}^m N_\Gamma(R-r;A,\gamma_i V)\le m
N_\Gamma(R;A,B)+c\,,\ee
 which yields the assertion. The claim for
weakly cocompact $\Gamma$ follows with Corollary~\ref{orbupcount}
and Lemma~\ref{orblowcount}.\qed\\

\section{Growth of conjugacy classes}

Let $M$ be a complete Riemannian manifold of nonpositive
sectional curvature with universal Riemannian covering manifold $\XX$, and
$\Gamma\subset\is(\XX)$ the group of deck transformations of the covering
projection $\XX\to M$.
It is well known that $\XX$ is a Hadamard manifold and
$\Gamma$ a discrete and torsion free group isomorphic to the fundamental group
of $M$. In this section, we
will derive a new asymptotic estimate for the growth rate of
geometrically distinct closed geodesics modulo free homotopy in
$M$. Notice that due to the occurrence of flat strips there can be
infinitely many closed geodesics in one free homotopy class.

We will only require that the group of deck transformations
$\Gamma$ of $M$ contains a hyperbolic axial element and possesses
infinitely many limit points.
Since we do not assume the manifolds to be compact or of finite
volume, we face certain difficulties which do not occur in the case of compact
manifolds treated  in \cite{K2}.
First, closed geodesics in $M$ may have arbitrarily small length.
Moreover, our weaker notion of hyperbolic geodesics includes the treatment of manifolds $M$ which are not necessarily of  geometric rank one, i.e. 
every geodesic in $\XX$ may bound a flat strip. In particular, we do not have a uniform upper bound on the width of hyperbolic axial isometries in $\Gamma$. 

For these two reasons  we encounter 
difficulties when trying to estimate the number of elements in $\Gamma$ which correspond to the same free homotopy class of closed geodesics. In particular, for lack of a uniform upper bound on the width of hyperbolic axial isometries, the argument in Lemma~5.4 of \cite{K2} 
cannot be directly adapted to our case. 

\begin{df} $\gamma, \gamma'\in\Gamma$ are said to be {\hd
equivalent} if and only if there exist $n,m\in \ZZ$ and
$\varphi\in\Gamma$ \st $(\gamma')^m=\varphi \gamma^n\varphi^{-1}$.
An element $\gamma_0\in\Gamma$ is called {\hd primitive} if it
cannot be written as a proper power $\gamma_0=\varphi^n$, where
$\varphi\in\Gamma$ and $n\ge 2$. \end{df} Each equivalence class
can be represented as
$$[\gamma]=\{\varphi \gamma_0^k\varphi^{-1}\;|\,
\gamma_0\in\Gamma,\;\gamma_0\
\mbox{primitive},\;k\in\ZZ,\;\varphi\in\Gamma\}\,.$$ It is easy to see that
the set of equivalence classes of axial elements in $\Gamma$ is in
one to one correspondence with the set of geometrically distinct
closed geodesics modulo free homotopy.
If $\gamma_0$ is a primitive axial isometry representing
$[\gamma]$, we have $$l([\gamma]):=\min \{l(\varphi)\;|\,
\varphi\in [\gamma]\}= l(\gamma_0)\,.$$ Then
$\ P(t):=\#\{[\gamma]\;|\, \gamma\in\Gamma\ \mbox{axial},\;
l([\gamma])\le t\}\ $ counts the number of closed geometrically
distinct geodesics of period $\le t$ modulo free homotopy, and 
$$P_{h}(t):=\#\{[\gamma]\;|\, \gamma\in\Gamma\
\mbox{hyperbolic axial},\; l([\gamma])\le t\}\le P(t)$$ 
the number of closed geometrically distinct hyperbolic geodesics of period $\le
t$ modulo free homotopy. If $W\subset\XX$ is an open set, we put
\be P(t;W)&:=& \#\{[\gamma]\;|\, \gamma\in\Gamma\
\mbox{axial},\;\Ax(\gamma)\cap W\ne \emptyset,\; l([\gamma])\le
t\}\,,\ \an\\
P_{h}(t;W)&:=& \#\{[\gamma]\;|\, \gamma\in\Gamma\
\mbox{hyperbolic axial},\;\Ax(\gamma)\cap W\ne \emptyset,\;
l([\gamma])\le t\}\,.\ee
For the remainder of this section we further fix a base point $\xo\in\XX$. Our
first lemma gives an easy upper bound for $P(t;W)$. 
\begin{lem}\label{upbound}
Let $W\subset\XX$ be a bounded open set. Then there exists a
constant $b>1$ \st $$P(t;W)\le b e^{\delta(\Gamma)t}\,.$$
\end{lem}
\prf\  Let $\gamma\in\Gamma$ be a primitive axial isometry, and
$x\in\Ax(\gamma)\cap W$. If $r:=\sup_{y\in W} d(\xo,y)$, then
$$d(\xo,\gamo)\le d(\xo,x)+d(x,\gamma x)+d(\gamma x,\gamma\xo)<
l(\gamma)+2r\,.$$ We conclude $P(t,W)\le N_\Gamma(t+2r)$  and, by
Corollary~\ref{orbupcount},\\[-1mm]

$\hspace{3cm} P(t;W)\le b e^{2r \delta(\Gamma)}e^{\delta(\Gamma)t}$.\qed\\

From here on we fix a hyperbolic axial isometry $h\in\Gamma$ and let 
$c_0:=c(h)\ge 0$ denote the width of $h$ (see Definition~\ref{hypaxiso}).  We further assume that $W\subset \XX$ is a bounded open set which contains the closure of a ball of radius $c_0$ centered at a point $y$ on an axis of $h$.

In order to bound $P_{h}(t;W)$ from below, we will need a few preliminary
lemmata. The first one gives a straightforward lower bound for $P_h(t;W)$. 
\begin{lem}\label{boundgivesopen}
There exist open  neighborhoods $\,U,V\subset\rand$ of $h^+$, $h^-$ with closures $\overline U$, $\overline V$ homeomorphic to closed balls and $\overline U\cap \overline V=\emptyset$ \st for all $t>0$ 
$$P_h(t;W)\ge \#\{[\gamma]\;|\, \gamma\in\Gamma \
\mbox{hyperbolic axial},\; \gamma^-\hspace{-0.5mm}\in U,\;\gamma^+\hspace{-0.5mm}\in V,\;l([\gamma])\le t\}\,.$$
\end{lem}
\prf\ Recall that $y$ is a point on an axis of $h$ and $W\supset \overline{B_y(c_0)}$. Let $\eps\in (0,1)$ be arbitrarily small with the property $B_y(c_0+\eps)\subseteq W$, and $U,V\subset\rand$ the corresponding neighborhoods of $h^-$, $h^+$ as in Lemma~\ref{hypgeod}. Then if $\gamma$ is hyperbolic axial with $l(\gamma)=l([\gamma])\le t$,  
$\gamma^-\in U$ and $\gamma^+\in V$, every axis $\sigma$ of
$\gamma$ satisfies $d(y,\sigma)< c_0+\eps$, hence
$\sigma\cap W\ne\emptyset$. We conclude that $[\gamma]$ is contained in the set 
$\ \{[\gamma]\;|\, \gamma\in\Gamma \
\mbox{hyperbolic axial},\; \Ax(\gamma)\cap W \ne\emptyset,\;l([\gamma])\le t\}\,.$\qed\\

The following lemma generalizes Lemma~5.4 in \cite{K2}. It gives the necessary  upper bound for the number of
$\gamma\in\Gamma$ which belong to the same equivalence class. 
\begin{lem}
Let $W\subset\XX$, $y\in\Ax(h)$, $c_0=c(h)$  and $U,V\subset\rand$ 
as in the previous lemma. Put 
$$\rho:=\frac14 \min_{\gamma\in\Gamma\setminus\{\idsm\}}\big(\inf_{x\in B_y(c_0+1)}  d(x,\gamma x)\big)\,.$$ 
Let $x\in\XX$,  and
$\gamma_0\in\Gamma$ a primitive hyperbolic axial element. 
Then there exists a constant $b>0$ depending only on $\rho$, $c_0$, $U$, $V$ and $N=\dim\XX$ \st for all $\,t>0$ 
\begin{eqnarray*}
\#\{\gamma=\varphi \gamma_0^k\varphi^{-1}\;|\,\varphi\in\Gamma,\;k\in\ZZ,\;
\varphi\gamma_0^-\in U,\,\varphi\gamma_0^+\in V,\,\ && \\
\varphi\Ax(\gamma_0)\cap B_x(\rho)\cap B_y(c_0+1)\ne\emptyset,\,\  l(\gamma)\le t \} & \le &  b \cdot t\,.
\end{eqnarray*}
\end{lem}
\prf\ Let $t_0:=l(\gamma_0)$ and $t\ge t_0$. Then for
$\gamma=\varphi \gamma_0^k\varphi^{-1}$ with $l(\gamma)\le t$ we
have $t\ge l(\gamma_0^k)=|k| t_0$, hence $|k|\le t/t_0$.

We remark that if $\varphi\gamma_0^-\in U$ and $\varphi\gamma_0^+\in V$, then
every axis $\sigma$ of $\gamma_0$ satisfies $d(y, \varphi\sigma)<c_0+1$ by  choice of $U$ and $V$. 

If $k\in\ZZ\setminus\{0\}$ is fixed, then
$\varphi\gamma_0^k\varphi^{-1}\ne \beta\gamma_0^k\beta^{-1}$
implies that $\beta^{-1}\varphi$ does not belong to the
centralizer of $\gamma_0$, in particular
$\varphi\Ax(\gamma_0)\neq\beta\Ax(\gamma_0)$.

Let $F_0\subset\Ax(\gamma_0)$ be a fundamental domain for the action of
$\langle \gamma_0\rangle$ on $\Ax(\gamma_0)$. If 
\begin{eqnarray*}
&&\varphi\gamma_0^-\in U,\,\  \varphi\gamma_0^+\in V,\,\ \varphi\Ax(\gamma_0)\cap
B_x(\rho)\cap B_y(c_0+1)\ne\emptyset,\,\ \an \\
&&\beta\gamma_0^-\in U,\,\ \beta\gamma_0^+\in V,\,\ \beta\Ax(\gamma_0)\cap
B_x(\rho)\cap B_y(c_0+1)\ne\emptyset\,,
\end{eqnarray*}
there exist $p,q\in F_0$ and $n,m\in\ZZ$ \st
$\varphi\gamma_0^n p$, $\beta\gamma_0^m q\in B_x(\rho)\cap B_y(c_0+1)$. Further\-more $g:=\varphi\gamma_0^n\ne
\beta\gamma_0^m=:f$ implies $d(p,q)\ge 2\rho(W)$ since \be 2
\rho(W)&\ge & d(gp,fq)\ge
d(gp,fp)-d(fp,fq)\\
&=& d(gf^{-1}gp,gp)-d(p,q)\ge 4 \rho(W)-d(p,q)\,.\ee 
If $d=\dim\Ax(\gamma_0)$, then $\vol(F_0)= (2c_0 +2)^{d-1}\cdot t_0\le 2^N (c_0+1)^{N-1}\cdot t_0\ $ and
$$ \vol(B_p(\rho)\cap \Ax(\gamma_0))=\omega_d\cdot\rho^{d}\qquad\forall\, p\in\Ax(\gamma_0)\,.$$ 
Put $\omega:=\min\{ \omega_d\;|\, 1\le d\le N\}$ and notice that $\rho^d\ge \min\{1,\rho^N\}$. Since the balls of radius $\rho$ centered at points in $F_0$ corresponding to different elements $\varphi\in\Gamma$ are disjoint, there are at most 
$$\frac{\vol(F_0)}{\omega_d\cdot \rho^d}\le \frac{2^N(c_0+1)^{N-1}\cdot t_0}{\omega\cdot\min\{1,\rho\}^{N}}$$ different
elements of the form $\varphi\gamma_0^k\varphi^{-1}$ \st
$\varphi\Ax(\gamma_0)\cap B_x(\rho)\cap B_y(c_0+1)\ne\emptyset$. 

The assertion now  follows from $\#\{k\in\ZZ\;|\,|k|\le t/t_0\}\le
2t/t_0$. \qed
\begin{cor}\label{numinequiclass}
There exists a
constant $e>0$ depending only on $c_0$, $U$, $V$ and $N=\dim\XX$ \st for all $\,t>0$ \begin{eqnarray*}
\#\{\gamma=\varphi \gamma_0^k\varphi^{-1}\;|\,\varphi\in\Gamma,\;k\in\ZZ,\;
\varphi\gamma_0^-\in U,\,\varphi\gamma_0^+\in V,\,\ && \\
 l(\gamma)\le t \} & \le &  e \cdot t\,.
\end{eqnarray*}
\end{cor}
\prf\ We use the notations from the previous  lemma and notice that by choice
of $U$, $V$ the conditions 
$\varphi\gamma_0^-\in U$ and $\varphi\gamma_0^+\in V$ imply that every axis
$\sigma$ of $\gamma_0$ satisfies  
\mbox{$\varphi\sigma\cap B_y(c_0+1)\ne \emptyset$.}

Since $\overline{B_y(c_0+1)}\subset\XX$ is compact, there exist finitely many balls $B_{x_i}(\rho)$, $1\le i\le m$, \st 
$$ \overline{B_y(c_0+1)}\subset\bigcup_{i=1}^m B_{x_i}(\rho)\,.$$
Hence if $\varphi\gamma_0^-\in U$ and $\varphi\gamma_0^+\in V$, there exists
$j\in \{1,2,\ldots, m\}$ \st\\ 
$\varphi\Ax(\gamma_0)\cap\big(B_{x_j}(\rho)\cap
B_y(c_0+1)\big)\ne \emptyset$. We conclude\\[-1mm]

$ \#\{\gamma=\varphi
\gamma_0^k\varphi^{-1}\;|\,\varphi\in\Gamma,\;k\in\ZZ,\;
\varphi\gamma_0^-\in U,\,\varphi\gamma_0^+\in V,\, l(\gamma)\le t
\}\le m\cdot a\cdot t\,.$\qed\\

We will now state two more lemmata in order to relate $P_{h}(t;W)\,$ to
$\,N_\Gamma(R;A,B)$ for appropriate sets $A,B\subset\rand$. 

\begin{lem}\label{nonintersect}
Let $\eps>0$ and $U,V\subseteq\rand$ be the corresponding disjoint
neighborhoods of $h^+$, $h^-$ as in Lemma~\ref{hypgeod}. Then there exist $\alpha>0$ and $R>0$
\st every $\gamma\in\Gamma$ with $d(\xo,\gamo)\ge R$,
$\angle_\xo(\gamo,h^+)<\alpha/2$ and
$\angle_\xo(\gamma^{-1}\xo,h^-)<\alpha/2$ satisfies $\gamma U\cap
V=\emptyset$.
\end{lem}
\prf\ Let $y\in\Ax(h)$, put $c:=c(h)+\eps$ and choose
$\delta>0$ \st $C_{y,h^+}^\delta\subset U$ and $C_{y,h^-}^\delta\subset
V$. From Lemma~2.5 and Lemma~2.8 in \cite{EbN} it follows that for any $x\in\XX$ there exist $T(x)>0$ and $\alpha(x)>0$ \st 
$$ C_{x,h^+}^{\alpha(x)}\setminus B_x(T(x))\subseteq C_{y,h^+}^{\delta/2}\qquad\mbox{and}\quad C_{x,h^-}^{\alpha(x)}\setminus B_x(T(x))\subseteq C_{y,h^-}^{\delta/2}\,.$$ 
Since $T$ and $\alpha$ depend continuously on $x$, for $T:=\max \{T(x)|\, x\in\{\xo\}\cup \overline{B_y(c)}\}$ and $\alpha:=\min \{\alpha(x)|\, x\in \{\xo\}\cup\overline{B_y(c)}\}$ we have
$$ C_{\xo,h^+}^{\alpha}\setminus B_\xo(T)\subseteq C_{y,h^+}^{\delta/2} \quad\mbox{and}\quad C_{x,h^-}^{\alpha}\setminus B_x(T)\subseteq C_{y,h^-}^{\delta/2}\quad\forall\, x\in \overline{B_y(c)}\cup \{\xo\}\,.$$

Moreover, since for $\gamma\in\Gamma$ and $x\in \overline{B_y(c)}$ we have $d(\gamma x,\gamma \xo)=d(x, \xo)\le c + d(y, \xo)\,$ and $d(\gamma^{-1} y,\gamma^{-1}\xo)= d(y,\xo)$, there exists $\ R>T+d(y,\xo)+c\,$ \st every $\gamma\in\Gamma$ with 
$d(\xo,\gamma\xo)>R$, $\angle_\xo(\gamma\xo,h^+)<\alpha/2$ and $\angle_\xo(\gamma^{-1}\xo,h^-)<\alpha/2$ satisfies 
$$\angle_y(\gamma x, h^+)<\delta/2\,\ \an\ \,
\angle_x(\gamma^{-1}y, h^-)<\delta/2\quad\ \forall\, x\in \overline{B_y(c)}\,.$$Now for $\xi\in U$ arbitrary 
there exists a hyperbolic geodesic $\sigma$ joining $h^-$ to $\xi$
with $d(y,\sigma)\le c$. If $x\in \overline{B_y(c)}$ is the orthogonal
projection of $y$ to $\sigma$, then
$\angle_x(\gamma^{-1}y,\xi)=\pi-\angle_x(\gamma^{-1}y,h^-)$. Considering the triangle with vertices $\gamma^{-1}y$, $x$ and $\xi$, we further have 
$\angle_{\gamma^{-1}y}(\xi,x)+\angle_{x}(\gamma^{-1}y,\xi)\le \pi\ $. We conclude
$$\angle_y(\gamma\xi,\gamma x)=\angle_{\gamma^{-1}y}(\xi,x)\le
\angle_x(\gamma^{-1}y,h^-)<\delta/2 \,, $$ 
and therefore
$\angle_y(\gamma\xi,h^+)\le \angle_y(\gamma\xi,\gamma
x)+\angle_y(\gamma x,h^+)<\delta$. In particular $\gamma\xi\in U$,
which proves $\gamma U\cap V\subseteq U\cap V=\emptyset$.\qed\\

The following lemma is due to G.~Knieper (\cite[Lemma~2.6]{K1}).
\begin{lem}\label{hypaxele}
Let $U,V\subset\rand$
be neighborhoods of $h^+$, $h^-$ with closures $\overline U$, $\overline V$ homeomorphic to closed balls and $\overline U\cap \overline
V=\emptyset$.  Then there exists a constant $\tau>0$, $n\in\NN$ \st
for all $\gamma\in \Gamma$ with $d(o,\gamo)<t$ and $\gamma U\cap
V=\emptyset$, the isometry $h^n\gamma h^n$ possesses an axis with
extremities in $U$ and $V$ and $l(h^n\gamma h^n)\le t+\tau$.
\end{lem}
\prf\ Let $n\in\NN$ \st $h^n(\ganz\setminus V)\subset U$ and
$h^{-n}(\ganz\setminus U)\subset V$. Since $\gamma U\cap
V=\emptyset$, we have $h^n\gamma h^n(\overline U)\subseteq h^n
\gamma U\subseteq h^n(\ganz\setminus V)\subset U$ and
$h^{-n}\gamma^{-1} h^{-n}(\overline V)\subseteq h^{-n}\gamma^{-1}
V\subseteq h^{-n}(\ganz\setminus U)\subset V$. Since $\overline U$ and $\overline V$ are each homeomorphic to a closed ball, Brouwer's
fixed point theorem implies that $h^n\gamma h^n$ assumes its
fixed points in $U$ and $V$.  Furthermore,

$d(o,h^n\gamma h^n\xo)\le
d(h^{-n}\xo,\xo)+d(o,\gamo)+d(\gamo,\gamma h^n\xo)\le 2 d(\xo, h^n\xo)+t$.\qed

\begin{thr}\label{main}
Let $M$ be a complete Riemannian manifold of nonpositive sectional curvature
with universal Riemannian covering manifold $\XX$, and $\Gamma\subset\is(\XX)$
the group of deck transformations of the covering projection.
Suppose $\Gamma$ contains a hyperbolic axial isometry $h$ and does not globally fix a point in $\rand$. Let $W\subset\XX$ be a bounded open set which
contains a closed ball of radius $c(h)\ge 0$ (as in Definition~\ref{hypaxiso}) centered at a point on an axis of $h$. Then
$$\delta(\Gamma)=\lim_{t\to\infty}\big(\frac1{t}\log P(t;W)\big)=
\lim_{t\to\infty}\big(\frac1{t}\log P_{h}(t;W)\big)\,.$$ 
For weakly cocompact $\Gamma$ there exist $b>0$ and $R>0$ \st for $t>R$
$$ \frac{1}{bt} e^{\delta(\Gamma)t}\le P_{h}(t;W)\le P(t;W)\le b e^{\delta(\Gamma)t}\,.$$
\end{thr}
\prf\ We choose $y\in\Ax(h)$ \st $W\supset \overline{B_y(c(h))}$ and open subsets $U,V\subseteq\rand$ of $h^+$, $h^-$  with closures $\overline U$, $\overline V$ homeomorphic to closed balls and $\overline U\cap
\overline V=\emptyset$ as in
Lemma~\ref{boundgivesopen}.  By Lemma~\ref{hypaxele}, there exists a
constant $\tau>0$ \st for any $t>0$ 
\be && \#\{\gamma\in\Gamma\;|\,
d(o,\gamo)\le t,\ \gamma U\cap V=\emptyset\}\\
&& \ \ \le\, \#\{\gamma\in\Gamma\;|\, \gamma\ \mbox{hyperbolic
axial with}\ \gamma^-\in U,\; \gamma^+\in V,\;l(\gamma)\le
t+\tau \}\,.\ee 
By Lemma~\ref{nonintersect}, there exist sets
$A\subset U$, $B\subset V$ and a constant $R>0$ \st for all $t>R$
$$\#\{\gamma\in\Gamma\;|\, d(o,\gamo)\le t,\,
\gamma U\cap V=\emptyset\}\ge N_\Gamma(t;A,B)-N_\Gamma(R;A,B)\,.$$
Using Lemma~\ref{boundgivesopen} and Corollary~\ref{numinequiclass} we conclude that for  $t>\tau+R$ 
\be P_{h}(t;W)&\ge & \#\{[\gamma]\;|\, \gamma\in\Gamma \
\mbox{hyperbolic axial},\; \gamma^-\in U,\;\gamma^+\in V,\;l([\gamma])\le t\}\\
&\ge & \frac{1}{e\cdot t}\#\{\gamma\in\Gamma\;|\, \gamma\ \mbox{hyperbolic
axial},\; \gamma^-\in U,\;\gamma^+\in V,\;l( \gamma)\le t\}\\
&\ge & \frac1{e\cdot t}\big( N_\Gamma(t-\tau;A,B) - b
e^{\delta(\Gamma)R}\big)\,,
\ee 
which, together with Corollary~\ref{twosets} and Lemma~\ref{upbound}, proves the assertion.\qed\\

We remark that for compact manifolds $M$  we may
choose a bounded  open set $W\subset\XX$ which contains a fundamental domain
for the action of $\Gamma$. Hence Theorem~\ref{main} implies
Theorem~B of G.~Knieper (\cite{K2}).

\vspace{1.0cm}
Gabriele Link\\
Mathematisches Institut II\\
Universit{\"a}t Karlsruhe\\
Englerstr.~2\\
76 128 Karlsruhe\\
e-mail:\ gabriele.link@math.uni-karlsruhe.de
\end{document}

%% file: commands.tex
\newcommand{\Zt}{\rm}

\newcommand{\ba}{\begin{array}}
\newcommand{\ea}{\end{array}}
\newcommand{\pot}{{\cal P}}
\newcommand{\curv}{\cal C}
\newcommand{\ddt} {\mbox{$\frac{\partial  }{\partial t}$}}
\newcommand{\hl}{\sf}
\newcommand{\hd}{\sf}

\newcommand{\Ad}{\mbox{\rm Ad}}
\newcommand{\Adsm}{\mbox{{\rm \scriptsize Ad}}}
\newcommand{\ad}{\mbox{\rm ad}}
\newcommand{\adsm}{\mbox{{\rm \scriptsize ad}}}
\newcommand{\diag}{\mbox{\rm Diag}}
\newcommand{\sect}{\mbox{\rm sec}}
\newcommand{\id}{\mbox{\rm id}}
\newcommand{\idsm}{\mbox{{\rm \scriptsize id}}}
\newcommand{\eps}{\varepsilon}

\newcommand{\aL}{\mathfrak{a}}
\newcommand{\bL}{\mathfrak{b}}
\newcommand{\mL}{\mathfrak{m}}
\newcommand{\kL}{\mathfrak{k}}
\newcommand{\gL}{\mathfrak{g}}
\newcommand{\nL}{\mathfrak{n}}
\newcommand{\hL}{\mathfrak{h}}
\newcommand{\pL}{\mathfrak{p}}
\newcommand{\uL}{\mathfrak{u}}
\newcommand{\lL}{\mathfrak{l}}

\newcommand{\kG}{{\tt k}}
\newcommand{\nG}{{\tt n}}

\newcommand{\Cart}{$G=K e^{\overline{\aL^+}} K$}
\newcommand{\Area}{\mbox{Area}}
\newcommand{\Hd}{\mbox{\rm Hd}}
\newcommand{\Hdim}{\mbox{\rm dim}_{\mbox{\rm \scriptsize Hd}}}
\newcommand{\Tr}{\mbox{\rm Tr}}
\newcommand{\bs}{{\cal B}}
\newcommand{\nc}{{\cal N}}
\newcommand{\MM}{{\cal M}}
\newcommand{\Ch}{{\cal C}}
\newcommand{\clCh}{\overline{\cal C}}
\newcommand{\Cnt}{\mbox{\rm C}}

\newcommand{\NN}{\mathbb{N}} \newcommand{\ZZ}{\mathbb{Z}}
\newcommand{\QQ}{\mathbb{Q}} \newcommand{\RR}{\mathbb{R}}
\newcommand{\KK}{\mathbb{K}} \newcommand{\FF}{\mathbb{F}}
\newcommand{\CC}{\mathbb{C}} \newcommand{\EE}{\mathbb{E}}
\newcommand{\XX}{X}
\newcommand{\HH}{I\hspace{-2mm}H}
\newcommand{\norm}{\Vert\hspace{-0.35mm}|}
\newcommand{\Sph}{\mathbb{S}}
\newcommand{\ganz}{\overline{\XX}}
\newcommand{\rand}{\partial\XX}
\newcommand{\prodrand}{\partial\XX_1\times\partial\XX_2} 
\newcommand{\regrand}{\partial\XX^{reg}}
\newcommand{\singrand}{\partial\XX^{sing}}
\newcommand{\Frand}{\partial\XX^F}
\newcommand{\Lim}{L_\Gamma}          
\newcommand{\cLim}{M_\Gamma}          
\newcommand{\Flim}{K_\Gamma}
\newcommand{\reglim}{L_\Gamma^{reg}}
\newcommand{\radlim}{L_\Gamma^{rad}}
\newcommand{\raylim}{L_\Gamma^{ray}}
\newcommand{\horinf}{\mbox{Vis}^{\infty}}
\newcommand{\horF}{\mbox{Vis}^B}
\newcommand{\Sml}{\mbox{Small}}
\newcommand{\SmlF}{\mbox{Small}^B}

\newcommand{\ifl}{\qquad\Longleftrightarrow\qquad}
\newcommand{\at}{\!\cdot\!}
\newcommand{\ging}{\gamma\in\Gamma}
\newcommand{\xo}{{o}}
\newcommand{\gamo}{{\gamma\xo}}
\newcommand{\gam}{\gamma}
\newcommand{\gax}{h}
\newcommand{\gxi}{{G\!\cdot\!\xi}}
\newcommand{\bd}{$(b,\Gamma\at\xi)$-densit}
\newcommand{\bt}{$(b,\theta)$-densit}
\newcommand{\cd}{$(\alpha,\Gamma\at\xi)$-density}
\newcommand{\be}{\begin{eqnarray*}}
\newcommand{\ee}{\end{eqnarray*}}

\newcommand{\an}{\ \mbox{and}\ }
\newcommand{\as}{\ \mbox{as}\ }
\newcommand{\diam}{\mbox{diam}}
\newcommand{\is}{\mbox{Isom}}
\newcommand{\Ax}{\mbox{Ax}}
\newcommand{\Fix}{\mbox{Fix}}
\newcommand{\Par}{F}
\newcommand{\Min}{\mbox{Fix}}
\newcommand{\vol}{\mbox{vol}}
\newcommand{\Td}{\mbox{Td}}
\newcommand{\piF}{\pi^B}
\newcommand{\piKM}{\pi^I}

\newcommand{\for}{\ \mbox{for}\ }
\newcommand{\pr}{\mbox{pr}}
\newcommand{\sh}{\mbox{sh}}
\newcommand{\shi}{\mbox{sh}^{\infty}}
\newcommand{\rank}{\mbox{rank}}
\newcommand{\supp}{\mbox{supp}}
\newcommand{\mass}{\mbox{mass}}
\newcommand{\kernel}{\mbox{kernel}}
\newcommand{\st}{\mbox{such}\ \mbox{that}\ }
\newcommand{\Stab}{\mbox{Stab}}
\newcommand{\Root}{\Sigma}
\newcommand{\Cone}{\mbox{C}}
\newcommand{\wrt}{\mbox{with}\ \mbox{respect}\ \mbox{to}\ }
\newcommand{\where}{\ \mbox{where}\ }

\newcommand{\con}{{\sc Consequence}\newline}
\newcommand{\rem}{{\sc Remark}\newline}
\newcommand{\prf}{{\sl Proof.\  }}
\newcommand{\qed}{$\hfill\Box$}

\newenvironment{rmk} {\newline{\sc Remark.\ }}{}  
\newenvironment{rmke} {{\sc Remark.\ }}{}  
\newenvironment{rmks} {{\sc Remarks.\ }}{}  
\newenvironment{nt} {{\sc Notation}}{}  

\newtheorem{satz}{\bf Theorem}

\newtheorem{df}{\sc Definition}[section]
\newtheorem{cor}[df]{\sc Corollary}
\newtheorem{thr}[df]{\bf Theorem}
\newtheorem{lem}[df]{\sc Lemma}
\newtheorem{prp}[df]{\sc Proposition}
\newtheorem{ex}{\sc Example}
